\documentclass{elsart}
\usepackage{amssymb,amsmath}
\begin{document}
\begin{frontmatter}

\title{Reduction of Rota's basis conjecture to a problem on
three bases}

\author{Timothy Y. Chow}
\address{Center for Communications Research, 805 Bunn Drive,
Princeton, NJ 08540}

\begin{abstract}
It is shown that Rota's basis conjecture
follows from a similar conjecture that involves
just three bases instead of $n$ bases.
\end{abstract}

\end{frontmatter}

\section{Introduction}
\label{intro}

In 1989, Rota formulated the following conjecture, which remains open.

\begin{conj}[Rota's basis conjecture]
\label{conj:rota}
Let $M$ be a matroid of rank~$n$ on $n^2$ elements
that is a disjoint union of $n$ bases
$B_1, B_2, \ldots, B_n$.
Then there exists an $n\times n$ grid~$G$
containing each element of~$M$ exactly once,
such that for every~$i$,
the elements of~$B_i$ appear in the $i$th row of~$G$,
and such that every column of~$G$ is a basis of~$M$.
\end{conj}

Partial results towards this conjecture may be found in
\cite{AB,Cha,Cho,Dr1,Dr2,GH,GW,HR,Pon,Wil,Zap}.
Now consider the following conjecture.

\begin{conj}
\label{conj:fixed}
Let $M$ be a matroid of rank~$n$ on $3n$ elements
that is a disjoint union of $3$ bases.
Let $I_1, I_2, \ldots, I_n$ be disjoint
independent sets of~$M$, with $0 \le |I_i| \le 3$ for all~$i$.
Then there exists an $n\times 3$ grid~$G$
containing each element of~$M$ exactly once,
such that for every~$i$,
the elements of~$I_i$ appear in the $i$th row of~$G$,
and such that every column of~$G$ is a basis of~$M$.
\end{conj}

The main purpose of the present note is to
make the following observation.

\begin{thm}
\label{thm:main}
Conjecture~\ref{conj:fixed} implies Conjecture~\ref{conj:rota}.
\end{thm}

Our proof is inspired by the proof of Theorem~4 in~\cite{Kei1}.

{\bf PROOF.}
Since Conjecture~\ref{conj:rota} is known if $n\le 2$,
we may assume that $n\ge 3$.
Let $M$ be given as in
the hypothesis of Conjecture~\ref{conj:rota}.
Define a {\em transversal\/} to be
a subset~$\tau \subseteq M$
that contains exactly one element from each~$B_i$.
Define a {\em double partition\/} of~$M$ to be
a pair $(\beta, \tau)$ where
$\beta = (\beta_1, \beta_2, \ldots, \beta_n)$
is a partition of~$M$ into $n$ pairwise disjoint bases~$\beta_i$
and $\tau = (\tau_1, \tau_2, \ldots, \tau_n)$
is a partition of~$M$ into $n$ pairwise disjoint transversals.
Given a double partition $(\beta,\tau)$, define
$$\mu(\beta,\tau) = \sum_{i\ne j} |\beta_i \cap \tau_j|.$$
Observe that if $\mu(\beta,\tau) = 0$ then
necessarily $\beta_i = \tau_i$ for all~$i$,
and then Rota's basis conjecture follows---just let
the $(i,j)$ entry of~$G$ be $B_i \cap \tau_j$.

So let $(\beta, \tau)$ be an arbitrary double partition
with $\mu(\beta,\tau) > 0$.
We show how to construct a double partition $(\beta', \tau')$
with $\mu(\beta', \tau') < \mu(\beta, \tau)$;
the proof is then complete, by infinite descent,
since by hypothesis there exists at least one double partition.
Since $\mu(\beta,\tau) > 0$,
there exist $\beta_i$ and $\tau_j$ with $i\ne j$ such that
$\beta_i \cap \tau_j \ne \emptyset$.
Since $n\ge 3$, there also exists~$k$
such that $i$, $j$, and~$k$ are all distinct.
It will simplify notation to assume that $i=1$, $j=2$, and $k=3$;
no generality is lost,
and it will be convenient to be able to reuse the index variables
$i$ and~$j$ below.
Let $S = \beta_1 \cup \beta_2 \cup \beta_3$,
let $T = \tau_1 \cup \tau_2 \cup \tau_3$,
and let $M' = M|S$ (i.e., $M$ restricted to the ground set~$S$).

For each $i$, let $I_i = B_i \cap T \cap S$.
Then $I_i$ is an independent subset of the matroid~$M'$,
and $|I_i| \le |B_i\cap T| \le 3$.
The $I_i$ are pairwise disjoint because the $B_i$ are pairwise disjoint.
Therefore we may apply Conjecture~\ref{conj:fixed}
to obtain an $n\times 3$ grid~$G'$ whose columns
$\beta'_1$, $\beta'_2$, and~$\beta'_3$ are disjoint bases of~$M'$
(and therefore are bases of~$M$)
and whose $i$th row contains the elements of~$I_i$.

To construct the desired double partition $(\beta', \tau')$,
let $\beta' = \beta$ except with
$\beta_1$, $\beta_2$, and~$\beta_3$
replaced with $\beta'_1$, $\beta'_2$, and $\beta'_3$ respectively.
Similarly, let $\tau' = \tau$ except with $\tau_1$, $\tau_2$, and~$\tau_3$
replaced with $\tau'_1$, $\tau'_2$, and~$\tau'_3$,
which are defined as follows.
Let $G''$ be any $n\times 3$ grid whose $i$th row
contains the elements of $B_i \cap T$ in some order,
and whose $(i,j)$ entry agrees with that of~$G'$
whenever that entry is in~$I_i$.
Clearly $G''$ exists (though it may not be unique).
Let $\tau'_j$ be the $j$th column of~$G''$, for $j=1,2,3$.

It is easily verified that
what we have done is to regroup the elements of~$M'$ into three new bases
and to regroup the elements of~$T$ into three new transversals in such a way
that the contribution to~$\mu(\beta', \tau')$ from intersections of
the new bases with the new transversals is reduced to zero,
and such that the total of the other contributions to~$\mu$ is unchanged.
Thus the overall value of~$\mu$ is reduced, as required.\qed

Careful inspection of the above proof shows that it is easily adapted
to prove a stronger statement than Theorem~\ref{thm:main}.
Let $C(k)$ denote the statement obtained by replacing `3' with `$k$'
throughout Conjecture~\ref{conj:fixed}.
Then the above argument, \emph{mutatis mutandis},
yields the following result.

\begin{thm}
\label{thm:stronger}
For any $\ell\ge k \ge 2$, $C(k)$ implies $C(\ell)$.
\end{thm}

In particular, proving $C(k)$ for any fixed~$k$
would prove Rota's basis conjecture
(in fact a stronger statement, namely $C(n)$)
for all $n$ greater than or equal to that fixed~$k$.

It is therefore natural to ask why we have
formulated Conjecture~\ref{conj:fixed} as $C(3)$ rather than as~$C(2)$.
The reason is that $C(2)$ is false.
The simplest counterexample is a well-known stumbling block
that is partly responsible for the fact that there is no known
general ``matroid union intersection theorem,''
i.e., a criterion for determining
the minimum number of common independent sets that
a set with two matroid structures on it
can be partitioned into.
Namely, take $M(K_4)$,
the graphic matroid of the complete graph on four vertices,
and let the $I_i$ be the three pairs of non-incident edges of~$K_4$.
Another counterexample arises from a matroid that Oxley~\cite{Oxl}
calls~$J$.  Representing~$J$ by vectors in Euclidean 4-space,
we can for example let
\begin{align*}
I_1 &= \{(-2, 3, 0, 1), (0, 0, 1, 1)\} \\
I_2 &= \{(0, 2, 0, 1),  (1, 0, 3, 1)\} \\
I_3 &= \{(1, 0, 0, 1), (0, 1, 2, 1)\} \\
I_4 &= \{(0, 1, 0, 1), (4, 0, 0, 1)\}
\end{align*}
It may be possible to construct other examples from
non-base-orderable matroids such as those in~\cite{Ing}.

Despite these counterexamples to~$C(2)$, we believe
that Conjecture~\ref{conj:fixed} is plausible.
Using a database of matroids with nine elements
kindly supplied by Gordon Royle~\cite{Roy},
we have computationally verified Conjecture~\ref{conj:fixed}
for the case $n=3$.

In an earlier version of this paper,
the formulation of Conjecture~\ref{conj:fixed} did not require
the $I_i$ to be independent.
A counterexample to that version
of the conjecture was found by Colin McDiarmid.
Take the complete graph on the
vertex set $\{1,2,3,4\}$,
and create an extra copy of the three edges incident to vertex~$4$.
Call the edges $12, 13, 14, 23, 24, 34, 14', 24', 34'$,
and let $I_1 = \{14, 14', 23\}$, $I_2 = \{24, 24', 13\}$,
and $I_3 = \{34, 34', 12\}$.
More generally, as pointed out by an anonymous referee,
if $k$ is odd, then a wheel with $k-1$ copies of each of its $k$ spokes
yields a counterexample to $C(k)$ if the $I_i$
are not required to be independent.

In closing, we speculate that Conjecture~\ref{conj:fixed}
might be provable using the following strategy.
First, develop a modified version of $C(2)$ that 
says that the conclusion holds provided certain
``obstructions'' (such as $M(K_4)$ and~$J$) are absent.
Then use Rado's theorem (12.2.2 of~\cite{Oxl}),
or a suitable strengthening of it,
to construct a first column of~$G$ in such a way
that the remaining $2n$ elements are obstruction-free.
Applying the modified version of~$C(2)$ would then
yield the desired result.
The analysis of obstructions should hopefully be tractable
since there are only 3 columns to consider.

\section{Acknowledgments}

I wish to thank Jonathan Farley, Patrick Brosnan, and James Oxley
for useful discussions, and a referee for correcting an error 
in my counterexample based on~$J$.

\end{document}